\documentclass[10pt]{article}

\usepackage{amsmath,amsthm}
\usepackage{makeidx}         
\usepackage{graphicx}        
\usepackage{multicol}        
\usepackage{amssymb}
\usepackage{latexsym}
\usepackage{url}

\theoremstyle{remark}

\theoremstyle{plain}

\newtheoremstyle{note}
  {3pt}
  {3pt}
  {}
  {}
  {\itshape}
  {:}
  {.5em}
  {}

\theoremstyle{note}

\newtheoremstyle{citing}
  {3pt}
  {3pt}
  {\itshape}
  {}
  {\bfseries}
  {.}
  {.5em}
  {\thmnote{#3}}

\theoremstyle{citing}

\newtheoremstyle{break}
  {9pt}
  {9pt}
  {\itshape}
  {}
  {\bfseries}
  {.}
  {\newline}
  {}

\theoremstyle{break}

\swapnumbers \theoremstyle{plain}

\let\lvert=|\let\rvert=|

\begin{document}

\begin{center}
{\LARGE\textbf{Empirical balanced truncation of nonlinear systems
\\}}\vspace {5mm}  \noindent

{\large \bf Marissa Condon$^{1,\dag}$} and {\large \bf Rossen I.
Ivanov$^{2,3,\ddag}$}\vskip 0.5cm  

\begin{tabular}{c}
\hskip-1cm $\phantom{R^R} ^{1}${\it School of Electronic
Engineering, Dublin City University,}
\\ {\it Glasnevin, Dublin 9, Ireland} \\
$\phantom{R^R} ^{2}${\it School of Mathematics, Trinity College
Dublin,}
\\ {\it Dublin 2, Ireland} \\
$\phantom{R^R}^{3}${\it Institute for Nuclear Research and Nuclear
Energy,}\\ {\it Bulgarian Academy of Sciences,} \\
{\it 72 Tzarigradsko chaussee, 1784 Sofia,
Bulgaria} \\
\\{\it $^\dag$e-mail: Marissa.Condon@dcu.ie}
\\ {\it $^\ddag$e-mail: ivanovr@tcd.ie}
\\
\hskip-.8cm
\end{tabular}
\vskip 0.5cm
\end{center}

\begin{abstract}
Novel constructions of empirical controllability and observability
gramians for nonlinear systems for subsequent use in a balanced
truncation style of model reduction are proposed.  The new gramians
are based on a generalisation of the fundamental solution for a
Linear Time-Varying system.  Relationships between the given
gramians for nonlinear systems and the standard gramians for both
Linear Time-Invariant and Linear Time-Varying systems are
established as well as relationships to prior constructions proposed
for empirical gramians. Application of the new gramians is
illustrated through a sample test-system.

{\bf AMS subject classification numbers}: 93B05, 93B07, 93B15, 93B40

{\bf Key Words}: controllability and observability gramians, model
reduction, balanced truncation, Lyapunov equation
\end{abstract}

\section{Introduction}

The development of effective model reduction techniques is of
paramount importance for all areas of engineering.  These include
control system design for nonlinear mechanical, chemical and
electronic engineering systems, the design of Radio-Frequency (RF)
integrated circuits and many others \cite{La} -- \cite{CI}.

In linear system theory (e.g. see \cite{ASG}, \cite{An} and the
references therein), the input-output interaction of a system is
characterized by the so-called \emph{gramian} matrices or
\emph{gramians}, which can be subsequently used in a model reduction
procedure, called balanced truncation \cite{ASG} -- \cite{ZDG}. For
general nonlinear systems the notion of gramians and balancing has
been derived from the more general concept of controllability and
observability (or \emph{energy}) functions \cite{Sch1} -- \cite{GS}.
However, the calculation of the energy functions is computationally
expensive and the result is rarely an explicit solution \cite{Ne},
\cite{Sch1} -- \cite{GS}. For these reasons, it is very difficult to
apply this method to large-scale problems \cite{La}. Several recent
research papers, \cite{La} followed by \cite{H1} -- \cite{H4}, have
presented a specific framework for the analysis and model reduction
of nonlinear models for the purpose of control termed
\emph{empirical balanced realization}. In the present paper, some
shortcomings of this approach as regards the determination of the
empirical gramians are detailed in Sections 3 and 4 and an improved
approach for the computation of the empirical gramians is suggested
in Section 5, \emph{Definitions 3,4}. Numerical tests are given in
Section 6.

\section{Empirical gramians and balanced truncation}

As in Lall et al. \cite{La}, the non-linear system under
consideration is of the form:
\begin{eqnarray}\label{NL1}
\dot{x}(t)&= &f((x(t),u(t)) \\\nonumber
 y(t)&=&h(x(t))
\end{eqnarray}
where $f:\mathbb{R}^n \times \mathbb{R}^p \rightarrow \mathbb{R}^n$
and $h: \mathbb{R}^n\rightarrow \mathbb{R}^q$ are nonlinear
functions, the function $u(t)\in \mathbb{R}^p$ is regarded as an
input signal to the system and the function $y(t)\in \mathbb{R}^q$
is an output signal. A simple idea, used extensively in the analysis
of autonomous nonlinear systems, is to compute a trajectory $x(t)$
on the time interval $[t_i,t_f]$ and to consider the integral
\cite{VD} $\int _{t_i}^{t_f} x(\tau)x(\tau)^T d\tau $ as an
approximation of the exact gramians for subsequent construction of
an appropriate projector (the superscript $T$ denotes
transposition). The method proposed in \cite{La} for general
nonautonomous systems stems from this basic idea. Data, taken either
from experiments or from numerical simulation and consisting of
sampled measurements of $x(t)$ and $y(t)$ , is used to parametrize
the trajectories for the nonlinear system.

The following constructions for empirical controllability and
observability gramians are then proposed in \cite{La}:

Let $\textbf{M}\equiv\{c_1,c_2,\ldots,c_s\}$   be a set of $s$
positive constants, $\textbf{T}^n \equiv\{T_1,T_2,\ldots,T_r\}$ --
be a set of $r$ orthogonal $n \times n$  matrices and $\textbf{E}^n
\equiv \{e_1,e_2,\ldots,e_n\}$  be the set of standard unit vectors
in $\mathbb{R}^n$ .

\emph{Definition 1}.  Let $\textbf{T}^p$, $\textbf{E}^p$  and
$\textbf{M}$  be given sets as described above. For the system
(\ref{NL1}) the empirical controllability gramian is defined as:

\begin{equation}
\hat{P}=\sum_{l=1}^{r} \sum_{m=1}^{s} \sum_{i=1}^{p}
\frac{1}{rsc^{2}_{m}} \int _0^{\infty} \Phi^{ilm} (t)dt \label{eq2}
\end{equation}
where $\Phi^{ilm} (t) \in \mathbb{R}^{n\times n}$  is given by
$\Phi^{ilm} (t)=(x^{ilm} (t)-\bar{x}^{ilm})(x^{ilm}
(t)-\bar{x}^{ilm})^T$ and $x^{ilm}(t)$ is the state of system
(\ref{NL1}) corresponding to the impulsive input $u(t)=c_m T_l e_i
\delta (t)$.  Here $\delta (t)$ denotes Dirac's delta function. The
mean $\bar{w}$ of a function $w \in L_1$ is given as:
\begin{equation}
\bar{w}= \lim _{t \rightarrow \infty} \frac{1}{t} \int _0^{t}
w(\tau)d\tau. \label{eq3}
\end{equation}

\emph{Definition 2}. Let $\textbf{T}^n$, $\textbf{E}^n$  and
$\textbf{M}$  be given sets as described above. For the system
(\ref{NL1}) the empirical observability gramian is defined as:
\begin{equation}
\hat{Q}=\sum_{l=1}^{r} \sum_{m=1}^{s} \frac{1}{rsc^{2}_{m}} \int
_0^{\infty} T_l \Psi^{lm} (t) T_l^T dt \label{eq4}
\end{equation}
where $\Psi^{lm} (t) \in \mathbb{R}^{n\times n}$  is given by
$\Psi^{lm}_{ij} (t)=(y^{ilm} (t)-\bar{y}^{ilm}
)^T(y^{jlm}(t)-\bar{y}^{jlm})$ and $y^{ilm}(t)$ is the output of
system (\ref{NL1}) corresponding to the initial condition $x^{ilm}
(0)=c_m T_l e_i$ with input $u=0$.

The purpose of using the sets $\textbf{M}$, $\textbf{T}^n$ and
$\textbf{E}^n$ in \emph{Definitions 1} and \emph{2} is an attempt to
ensure that the entire region of feasible values of initial
inputs/states is covered and probed. The set $\textbf{E}^n$ defines
the standard directions and the set $\textbf{T}^n$ defines
'rotations' of these directions. The set $\textbf{M}$ introduces
different scales for each direction of the initial states/inputs.

In what follows several shortcomings associated with
\emph{Definitions 1} and \emph{2} are brought to light and novel
proposals for improvement are suggested.

\section{Linear time-varying systems}

An examination of Linear Time-Varying Systems (LTVS) in the context
of model reduction is both nontrivial and instructive.  The
controllability gramian proposed in \emph{Definition 1} for a
non-autonomous system $\dot{x}=f(x,t)$ does not yield the standard
controllability gramian for such systems \cite{Mo}, \cite{VD}.
Furthermore, the derivation of the standard gramian for LTVS
provides a motivation for the new improved constructions suitable
for nonlinear systems. In what follows, for simplicity, only
one-dimensional inputs and outputs are considered, i.e. $p=q=1$   in
(\ref{NL1}), (\ref{eq2}). Consider a LTVS:
\begin{eqnarray}\label{eq5}
\dot{x}(t)&= &A(t)x(t)+B(t)u(t) \\\nonumber
 y(t)&=&C(t)x(t)
\end{eqnarray}
The fundamental solution of (\ref{eq5}) is defined as the solution
of:
\begin{equation}
\dot{\Theta}(t)=A(t)\Theta (t),\qquad  \Theta (0)=I \label{eq6}
\end{equation}
where $I$  is the corresponding identity matrix.  For example, if
$A$ is a constant matrix, (as for the linear time invariant system
-- LTIS) then one simply recovers the very well known solution
$\Theta (t)=\exp (At)$. The general solution of (\ref{eq5}) is:
\begin{equation}
x(t)=\Theta (t)  \left( \Theta ^{-1} (t_0)x(t_0)+\int _{t_0}
^{t}\Theta ^{-1}(\tau)B(\tau)u(\tau)d\tau \right) \label{eq7}
\end{equation}
Now let $t_0 \rightarrow -\infty$, $t=0$  and $x(-\infty)=0$.  From
(\ref{eq7}) it follows:
\begin{equation}
x(0)=\int _{-\infty} ^{0}\Theta ^{-1}(\tau)B(\tau)u(\tau)d\tau=\int
_{0} ^{\infty}\Theta ^{-1}(-\tau)B(-\tau)u(-\tau)d\tau \label{eq8}
\end{equation}
and as usual, one can define a Controllability operator:
\begin{equation}
\textbf{C}: L_2([0,\infty ))\rightarrow \mathbb{R}^n \qquad
\text{as} \qquad \int _{0} ^{\infty} d\tau \Theta
^{-1}(-\tau)B(-\tau)\bullet \label{eq9}
\end{equation}
and Controllability gramian as:
\begin{equation}
P= \int _{0} ^{\infty} \Theta ^{-1}(-\tau)B(-\tau) B^T (-\tau)\Theta
^{-1T}(-\tau) d\tau \label{eq10}
\end{equation}
From  (\ref{eq7}) with $t_0=0$  and $u\equiv 0$ it follows
$y(t)=C(t)\Theta (t)x(0)$ and therefore the Observability operator
can be defined as:
\begin{equation}
\textbf{O}: \mathbb{R}^n \rightarrow L_2([0,\infty )) \qquad
\text{as} \qquad \textbf{O}=C(t) \Theta (t) \label{eq11}
\end{equation}
and the Observability gramian is:
\begin{equation}
Q= \int _{0} ^{\infty} \Theta ^{T}(\tau)C^T(\tau) C (\tau)\Theta
(\tau) d\tau \label{eq12}
\end{equation}

Strictly speaking, the gramians for LTVS must depend on $t$  as
shown in \cite{Mo}, \cite{VD}.  However, for the purposes of model
reduction, constant gramians are preferred and the constant versions
(\ref{eq10}) and (\ref{eq12}) are used as approximations. The
expressions in (\ref{eq10}) and (\ref{eq12}) are generalisations of
the gramians for LTIS where $\Theta (t)=\exp (At)$ \cite{La}.

\section{Bilinear representation of nonlinear systems}

Another very interesting class of nonlinear systems that it is
instructive to examine are the bilinear systems; moreover a wide
class of nonlinear systems (subject to suitable restrictions--
\cite{Mo}, \cite{Ph}, \cite{CI}), may be represented in a bilinear
form. The bilinear system is also interesting because there is an
exact solution when the input is a delta-function and thus the
gramians (\ref{eq2}) and (\ref{eq4}) can be tested explicitly.
Consider the following bilinear system:
\begin{eqnarray}\label{eq13}
\dot{\hat{x}}(t)&= &\hat{A}(t)\hat{x}(t)+\hat{N} \hat{x}(t)u(t)+
\hat{B}u(t)
\\\nonumber
y(t)&=&\hat{C}\hat{x}(t)
\end{eqnarray}

Again, it is assumed that all the eigenvalues of $\hat{A}$  have
negative real parts.  Let the sets employed in \emph{Definition 1}
be as follows: $\textbf{M}\equiv\{c_1,c_2,\ldots,c_s\}$,
$\textbf{T}\equiv\{1\}$  and $\textbf{E}\equiv\{1\}$ since $p=q=1$.
Thus the inputs to the system are of the form $u_0(t)=c_m \delta
(t)$ . The solution to (\ref{eq13}) with an input $u_0(t)=c_m \delta
(t)$ is:
\begin{equation}
\hat{x}^{11m}(t)= e^{\hat{A}t}\left(
I+\frac{c_m}{2}\hat{N}+\frac{c_m ^2}{4}\hat{N}^2+\ldots
\right)\hat{B}c_m\theta (t) \label{eq14}
\end{equation}
where $\theta (t)$ is the unit step function. Note that the sum in
(\ref{eq14}) is finite since $\hat{N}$ is nilpotent by construction
\cite{Mo}, \cite{Ph}.  ($\hat{x}^{11m}(t)$ corresponds to
$\hat{x}^{ilm}(t)$ with $i=1$, $l=1$).  Following from
\emph{Definition 1}, the Controllability gramian is therefore:
\begin{equation}
P_{BL}= \int _{0} ^{\infty} e ^{\hat{A}\tau}\bar{B}_N \bar{B}^T _N e
^{\hat{A}^T \tau} d\tau \label{eq15}
\end{equation}
where
\begin{equation}
\bar{B}_N \bar{B}^T _N= \sum_{m=1} ^{s} \left(
I+\frac{c_m}{2}\hat{N}+\frac{c_m ^2}{4}\hat{N}^2+\ldots
\right)\hat{B} \hat{B}^T \left( I+\frac{c_m}{2}\hat{N}+\frac{c_m
^2}{4}\hat{N}^2+\ldots \right)^T .\label{eq16}
\end{equation}
Since the bilinear system (\ref{eq13})  assumes a linear form when
the input is zero, the Observability gramian is as usual:
 \begin{equation}
Q_{BL}= \int _{0} ^{\infty} e ^{\hat{A}^T \tau}\hat{C}^T \hat{C} e
^{\hat{A}\tau} d\tau \label{eq17}
\end{equation}
It is not difficult to prove that the gramians in (\ref{eq15}) and
(\ref{eq17}) are solutions to the following Lyapunov Equations:
\begin{eqnarray}\label{eq18}
\hat{A} P_{BL}+P_{BL}\hat{A}^T+\bar{B}_N \bar{B}^T _N &= &0
\\\nonumber
\hat{A}^T Q_{BL}+Q_{BL}\hat{A}+\hat{C}^T \hat{C}&=& 0
\end{eqnarray}
However, there are the following problems with the gramians in
(\ref{eq15}) and (\ref{eq17}).  Firstly, they do not relate to the
known gramians for the bilinear systems \cite{AlB} -- \cite{DA},
\cite{CI}.   Secondly, (\ref{eq14}) suggests that the Krylov space
for the Controllability operator is
$\text{span}\{\hat{A}^{p_1}\hat{N}^{p_2}\hat{B}\}$ for $p_i\geq 0$.
However, the known Krylov space \cite{Ph} is $\text{span}\{\hat{B};
\hat{A}^{p_1}\hat{B};\hat{A}^{p_1}\hat{N}
\hat{A}^{p_2}\hat{B};\ldots; \hat{A}^{p_1}\hat{N}
\hat{A}^{p_2}\hat{N} \ldots\hat{A}^{p_k}\hat{B}\}$ for $p_i
> 0$.

\section{Nonlinear systems}

The nonlinear system in (\ref{NL1}) has a rather general form.  In
\cite{H1},  \cite{H3} it is suggested that the use of the empirical
gramians (\ref{eq2}) and (\ref{eq4}) is limited only to
control-affine systems. Indeed, for example, for a system, depending
quadratically on the input, the square of the Delta-function cannot
be defined.

For the present analysis, let the nonlinear systems be of the form:
\begin{eqnarray}\label{eq19}
\dot{x}(t)&= &f(t, x(t))+ B(t)u(t) \\\nonumber
 y(t)&=&h(t,x(t))
\end{eqnarray}
It contains two terms: a \emph{dynamical} term (or drift term)
$f(t,x(t))$ and a \emph{source} term (or diffusion term) $B(t)u(t)$
. Clearly, LTVS systems are of the form in (\ref{eq19}).

Instead of considering different inputs and 'mean values' as in
\emph{Definitions 1} and \emph{2}, it is more natural to analyse the
system in a vicinity of an equilibrium point when $u(t)=0$ .
Consider the vicinity of an isolated asymptotically stable
equilibrium point (steady--state solution) which is supposed to be a
constant solution and is chosen for simplicity at $x=0$, i.e.
$f(t,0)\equiv 0$. It is also assumed that the system does not leave
the region of attraction of this equilibrium point when the input is
applied for the initial data used.  If the system exhibits multiple
steady--state solutions, then the analysis may be applied separately
in the vicinity of each solution provided that extra care is taken
to ensure that the system does not leave the region of attraction of
the corresponding (asymptotically stable) equilibrium point.  Of
course, the constructed gramians will therefore only provide a basis
for reduction locally in the vicinity of the selected equilibrium
point.

In this work, it is proposed to make use of an approximation for the
most natural object -- the fundamental solution $\Theta$ of
(\ref{eq19}) that would generalize the $\exp(At)$  term for linear
systems. This is reasonable since the projection Krylov spaces for
linear systems are generated by their fundamental solution
$\exp(At)$.  The constructions would, in general, depend on $\Theta$
for negative times which is unavoidable.  For linear systems, of
course, there is a simplification since $\left(
e^{A(-t)}\right)^{-1}\equiv e^{At}$ so this does not present a
limitation but in general, $\Theta ^{-1} (-t) \neq \Theta (t)$  ,
cf. (\ref{eq10}).

Let $x^{ilm}(t)$ be the solution of (\ref{eq19}) with  $u\equiv 0$:
\begin{equation}
\dot{x}(t)= f(t, x(t)) \label{eq20}
\end{equation}
and with initial condition:
\begin{equation}
x^{ilm}(0) = c_m T_l e_i \label{eq21}
\end{equation}
It is assumed that the initial condition (\ref{eq21}) does not take
the system outside the region of attraction of the equilibrium point
$x=0$.  Then the 'state-space average' of the 'nonlinear'
fundamental solution may be defined as:
\begin{equation}
\langle \Theta (t) \rangle = \frac{1}{rs} \sum_{m=1}^{s}
\sum_{l=1}^{r} \sum_{i=1}^{n} \frac{1}{c_m}x^{ilm}(t)e_i ^T T^T _l
\label{eq22}
\end{equation}
where the sets $\textbf{M}$, $\textbf{T}^n$, $\textbf{E}^n$
previously defined for \emph{Definitions 1} and \emph{2} are
employed. Indeed, for a LTVS, $x^{ilm}(t) = \Theta (t) c_m T_l e_i$
and therefore $\langle \Theta (t) \rangle \equiv \Theta (t)$.

The following constructions of empirical controllability and
observability gramians for the nonlinear system (\ref{eq19}) are now
suggested:

\emph{Definition 3}. For the system in (\ref{eq19}), the
\emph{nonlinear} controllability gramian is defined as:
\begin{equation}
\tilde{P}= \int _{0} ^{\infty}\langle \Theta (-\tau) \rangle
^{-1}B(-\tau) B^T (-\tau) \langle\Theta (-\tau) \rangle ^{-1T} d\tau
\label{eq23}
\end{equation}
where $\langle \Theta (t) \rangle$  is as described in (\ref{eq22}).

Of course, this construction requires that $\langle \Theta (-\tau)
\rangle$  is invertible for all $\tau \geq 0$. (\ref{eq23}) is
obviously a generalisation of (\ref{eq10}).

\emph{Definition 4}.  For the system in (\ref{eq19}) the
\emph{nonlinear} observability gramian is defined as:
\begin{equation}
\tilde{Q}= \int _{0} ^{\infty}z^T (\tau) z (\tau) d\tau \label{eq24}
\end{equation}
where $z(\tau)\in \mathbb{R}^n$ is given by:
\begin{equation}
z(t)= \frac{1}{rs} \sum_{i,l,m} \frac{1}{c_m}y^{ilm}(t)e_i ^T T^T _l
\nonumber
\end{equation}
and $y^{ilm}(t)$  is the output which corresponds to an initial
state $x^{ilm}(0) = c_m T_l e_i$  and a zero source term.  The
motivation for this construction is as follows:

For a linear output $y(t)=C(t)x(t)$, since $\langle \Theta (t)
\rangle \equiv \Theta (t)$ the observability gramian (\ref{eq12})
is:
\begin{equation}
Q= \int _{0} ^{\infty} \langle \Theta (\tau) \rangle ^{T} C^T(\tau)
C (\tau)\langle \Theta (\tau)\rangle d\tau \label{eq25}
\end{equation}
Since \begin{equation}C (\tau)\langle \Theta
(\tau)\rangle=\frac{1}{rs} \sum_{i,l,m}
\frac{1}{c_m}C(t)x^{ilm}(t)e_i ^T T^T _l= \frac{1}{rs} \sum_{i,l,m}
\frac{1}{c_m}y^{ilm}(t)e_i ^T T^T _l=z(t)\nonumber \end{equation}

\noindent then the construction in (\ref{eq24}) is confirmed as a
generalisation of (\ref{eq12}).

Both gramians (\ref{eq23}) and (\ref{eq24}) when applied to LTVS (or
LTIS) thus result in the usual gramians i.e. (\ref{eq10}) and
(\ref{eq12}). This confirms the motivation for their use in
preference to (\ref{eq2}) and (\ref{eq4}).

\section{Illustrative numerical example}

The circuit employed is the nonlinear $RC$ ladder shown in Fig. 1
(frequently employed as a test circuit for model reduction
techniques \cite{Ph} -- \cite{Do}, \cite{CI}).  The example enables
comparisons to be made between the existing formulations for
empirical gramians and those proposed in this contribution.  The
nonlinear resistors (a diode in parallel with a unit resistor) have
the constitutive relation $i(v)=(e^{40v}-1)+v$ (where $i$ represents
current and $v$ represents voltage). The capacitors have unit
capacitance.  The input is a current source $u(t)=e^{-t}$  entering
node 1 and the output is the voltage taken at node 1, Fig 2(a). This
is an example of a gradient system (e.g. according to the definition
in \cite{McL}), since the equations describing the system may be
written in the form:
\begin{eqnarray}\label{eq26}
\dot{v}&= &- \nabla V + B u(t) \\\nonumber
y&=&Cv\equiv v_1 (t)
\end{eqnarray}
where $B=[1\quad 0 \quad \ldots \quad 0]^T$, $C=B^T$ and
\begin{equation}
V(v)= \frac{1}{40} e^{40v_1} -v_1 + \frac{v^2 _1}{2} + \sum _{k=1}
^{n-1} \left ( \frac{1}{40} e^{40(v_k -v_{k+1})} -(v_k -v_{k+1})+
\frac{(v_k -v_{k+1})^2}{2} \right).\label{eq27}
\end{equation}

The function $V(v)$ represents a strong Lyapunov function for the
gradient system as described in \cite{McL}.  This then enables the
application of Lyapunov stability criterion to show that $v=0$  is
an asymptotically stable equilibrium point of the system (when the
source is set to zero).

The number of nodes in the system is  $n=30$.  The time interval
chosen for consideration is  $t\in [0, 1]$.  The reduction of the
original system to a system of order 3 is implemented using several
different methods.

In order to compare the new gramians with the existing constructions
for empirical gramians (\emph{Definitions 1,2}), a bilinear
representation \cite{Mo}, \cite{Ph}, \cite{CI} of the system in
(\ref{eq26}) -- (\ref{eq27}) is employed. The reason for doing this
is that an exact solution exists for a bilinear system when
subjected to impulsive inputs.  This is of importance in the
formation of the gramian as specified in \emph{Definition 1} as it
necessitates subjecting the system to impulsive inputs.  A bilinear
approximation with two terms in the Taylor's series expansion is
employed.  The resultant bilinear model is of order $30+30^2=930$.
For information, the Root Mean Square (RMS) error between the result
from the  nonlinear model (\ref{eq26}) and the full order-- $930$
bilinear approximation (\ref{eq13}) is $1.0 \times 10^{-2}$, Fig
2(b).

As a benchmark, consider the simplest reduced model (of order 3)--
that which employs only the linear part of the bilinear
approximation to form the gramians necessary for balancing.  To be
specific, the gramians employed are the solutions of the following
Lyapunov equations:
\begin{equation}\label{eq28}
\hat{A} P_{BL}+P_{BL}\hat{A}^T+\hat{B}\hat{B}^T = 0, \qquad
\hat{A}^T Q_{BL}+Q_{BL}\hat{A}+\hat{C}^T \hat{C}= 0.
\end{equation}
The RMS error in comparison to the full order-- $930$ bilinear model
is $2.6 \times 10^{-2}$, Fig 2(c).

Now consider the use of the gramians (\ref{eq18}) formed on the
basis of \emph{Definitions 1} and \emph{2} with
$\text{dim}(\hat{x})=930 $, $\textbf{M} \equiv \{-5,-0.5,-1,-0.1,
0.1, 0.5, 1, 5\}$, $\textbf{T}^{930}=\{I\}$.  The RMS error in
comparison to the full order-- $930$ bilinear model is $7.5 \times
10^{-2}$, Fig 2(d). Moreover, it is observed that when $\textbf{M}
\equiv \{c_1\}$, i.e. consisting of only one constant, the reduction
process is ill-defined for some values of $c_1$, e.g. $c_1 =0.20$;
$0.22$; i.e. the reduced model is unstable.

Finally, consider the case where the gramians formed on the basis of
\emph{Definitions 3} and \emph{4} are employed for reduction
purposes. The integration over $\tau$  in these constructions is
replaced by a discrete summation.  The resulting RMS error (in
comparison to the original model) is $5.3\times 10^{-5}$, Fig 2(e).
This indicates the superiority of the novel constructions for the
purposes of model reduction via a balancing technique.

\section{Conclusions}

The paper has proposed new constructions for empirical gramians for
subsequent use in a method of model reduction based on 'balancing'.
The important new concept involved in the formation of the novel
empirical gramians, (\ref{eq23}) and (\ref{eq24}), is that of a
'state-space average' of the 'nonlinear' fundamental solution
(\ref{eq22}).

The method is successful if the state-space average of the nonlinear
fundamental solution is well defined.  Of course, this is not the
case for all nonlinear systems as the solution of (\ref{eq20}) may
not exist or may only exist for specific choices of the initial
data. However, the method is applicable for systems for which the
nonlinearities are not too severe, e.g. for the so-called 'weakly'
nonlinear systems as described in \cite{Ph}.  For such systems, it
is expected that the 'nonlinear' fundamental solution is 'close' to
the exponential form that corresponds to the fundamental solution
for a linear system. The new empirical gramians coincide with the
usual gramians for both LTVS and LTIS.


\newpage
\begin{figure}
\centering
\includegraphics[height=4cm]{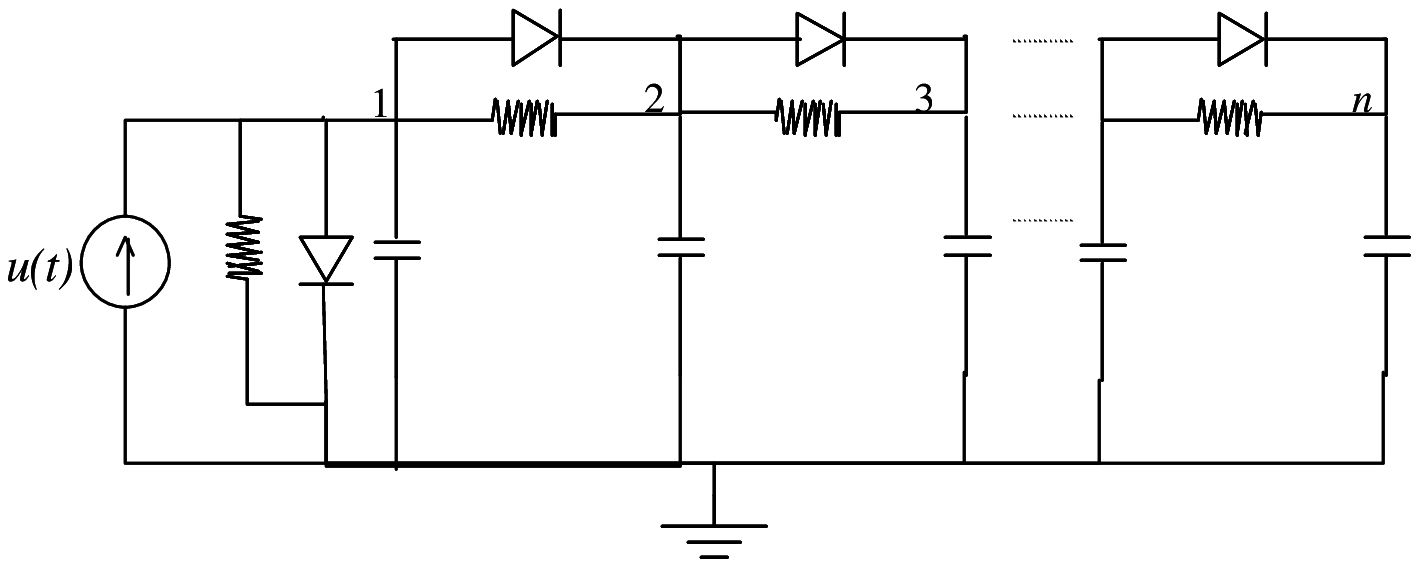}
%
%
\caption{Nonlinear circuit}
\label{fig:1}       
\end{figure}

\newpage

\begin{figure}
\centering
\includegraphics[height=9.5cm]{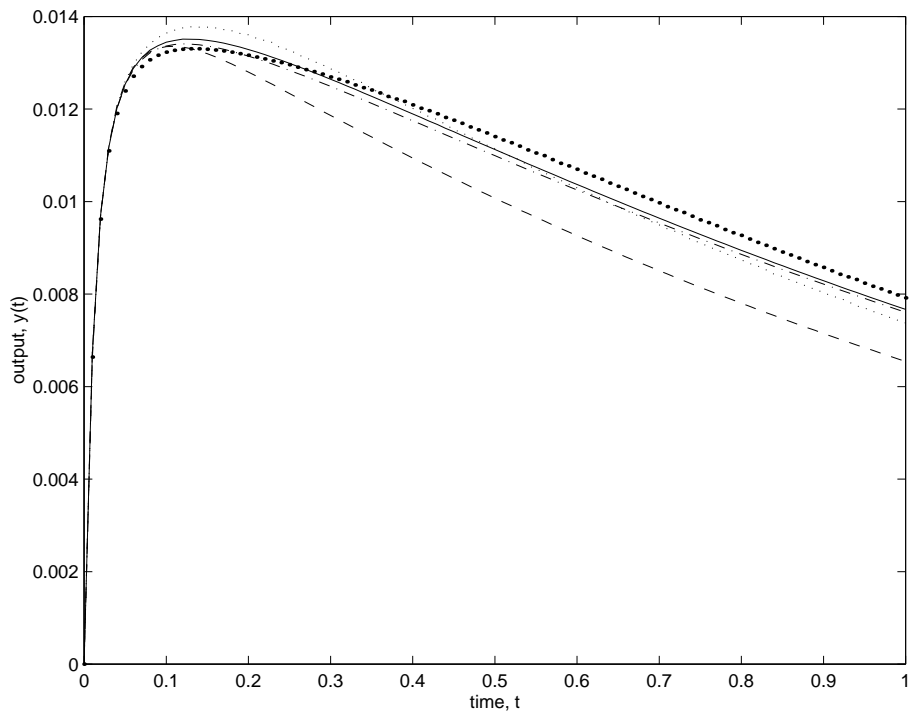}
%
%
\caption{Comparison between output from nonlinear model and
reduced-order models: (a) Solid line -- Nonlinear model
(\ref{eq26}); (b) Dash-dotted line -- Bilinear approximation
(\ref{eq13});  (c) Points -- Reduced bilinear system with gramians
based only on linear part of bilinear approximation (\ref{eq28});
(d) Dashed line -- Reduced-model with gramians (\ref{eq18}) based on
\emph{Definitions 1,2}. (e) Dotted line -- reduced-order model where
the reduction is based on the novel Empirical gramians --
(\ref{eq23}) and (\ref{eq24}). }
\label{fig:2}       
\end{figure}

\end{document}